\documentclass[12pt,a4paper]{article}
\usepackage{hyperref}
\usepackage{amssymb}
\usepackage{amsthm}
\usepackage{amsmath}
\usepackage{color}
\usepackage{graphicx}
\usepackage{comment}
\usepackage[margin=3cm]{geometry}
\usepackage{todonotes}
\usepackage{floatrow}

\newtheorem{theorem}{Theorem}[section]

\newtheorem{lemma}[theorem]{Lemma}
\newtheorem{proposition}[theorem]{Proposition}

\theoremstyle{definition}
\newtheorem{definition}[theorem]{Definition}
\newtheorem{assumption}[theorem]{Assumption}

\theoremstyle{remark}
\newtheorem{remark}[theorem]{Remark}

\newcommand{\C}[1]{\mathcal{C}_{#1}}
\newcommand{\G}[1]{\mathcal{G}_{#1}}
\newcommand{\pv}{\mathrm{p.v.}}
\newcommand{\sech}{\mathrm{sech}}
\newcommand{\csch}{\mathrm{csch}}

\newcommand{\ri}{\mathrm{i}}

\newcommand{\symb}{\mathrm{symb}}

\title{Variational Instability for Irrotational Water Waves in Finite Depth}

\author{Florian Kogelbauer}

\begin{document}


\date{\today}%

\maketitle

\begin{abstract}
We study the variational stability of solutions to the periodic irrotational gravity water-wave problem in finite depth with fixed conformal mean depth. We derive a finite-depth Plotnikov transformation for the second variation, including the zero-mode correction specific to finite depth. For fixed conformal mean depth, the trivial solution is strictly variationally stable for $\mu>\tanh(kh)/k$, becomes degenerate at $\mu=\tanh(kh)/k$, and is variationally unstable for smaller $\mu$. Along the primary non-trivial branch, a Lyapunov--Schmidt expansion and a direct evaluation of the Hessian on the bifurcating branch yield a negative direction for every finite depth and sufficiently small non-zero amplitude.
\end{abstract}

\section{Introduction}

Traveling wave solutions to the fully nonlinear water-wave problem exhibit several well-known instability mechanisms. A significant body of literature is dedicated to instability results for free-boundary flows. Most notably, for sufficiently deep water, in \cite{benjamin1967disintegration} it was shown that irrotational water waves in a small-amplitude regime exhibit the so-called
\textit{sideband instability} or \textit{Benjamin--Feir instability}. This phenomenon indicates that periodic waves are unstable under perturbations of a different period \cite{bridges1995proof}. For various classes of exact solutions, different kinds of instabilities have been rigorously proved.
For instance, instability has been established for Gerstner’s wave \cite{gerstner1809theorie,leblanc2004local},
equatorial waves \cite{constantin2013instability,genoud2014instability} and edge waves \cite{ionescu2014instability}. Similarly, transverse instability of Stokes waves in finite depth has been investigated  in \cite{creedon2024transverse} and classical numerical studies such as \cite{mclean1982instabilities} have provided significant insights into instability mechanisms. Furthermore, recent numerical analyses of normal-mode perturbations for two-dimensional periodic finite-amplitude gravity waves have been reported in \cite{francius2017two}.\\
In contrast, traveling wave solutions to model equations for water wave systems often display strong stability properties. Orbital stability has been proven for periodic wave solutions to equations of Korteweg–de Vries type \cite{mckean1977stability,pava2008positivity}, the Camassa-Holm equation \cite{lenells2005stability} and the Benjamin--Bona--Mahony equation \cite{angulo2011regularized}. Stability results for solitary traveling wave solutions have been obtained in \cite{ALBERT1987343,albert1992positivity} for long wave models, in \cite{albert1991total} for internal waves and in \cite{benjamin1972stability} for the Korteweg-de-Vries equation. 
General theorems on stability, applicable to a variety of systems, are summarized in  \cite{grillakis1987stability}.
These model equations often rely on asymptotic approximations, such as small-amplitude, shallow-water, or long-wave regimes, which can substantially affect  their stability properties.  For discussions on the linear stability of irrotational waves using the Hamiltonian formulation in the velocity potential framework, see \cite{mackay1986stability,zakharov1968stability}.
Surface tension effects are neglected in this study. However, when surface tension is considered, stability and energetic stability can be shown for capillary-gravity waves in both periodic \cite{buffoni2005minimization} and solitary waves \cite{buffoni2004existence,mielke2002energetic}. Moreover, stability results for flows with resistance to surface stretching and bending, as in \cite{buffoni2013stability}, reinforce the intuition that surface tension has a stabilizing influence. We also note the stability analysis for waves with non-zero vorticity, as detailed in \cite{constantin2007stability}.\\
In this work, we investigate the variational properties of periodic traveling waves for the irrotational free-boundary Euler equations under gravity in finite depth. More precisely, we study the  definiteness of the second variation of a variational principle under perturbations of the free-surface variable while the conformal mean depth is fixed. We call the presence of a negative direction \textit{variational instability}. This terminology refers only to the constrained steady variational problem considered here and no implication for spectral, linear, or nonlinear stability of the time-dependent water-wave problem is asserted.\\
We rely on a reformulation of the governing equations as a one-dimensional pseudo-differential equation coupled to a scalar constraint, obtained in \cite{constantin2016global,constantin2011steady} for traveling waves with constant vorticity. The pseudo-differential equation is the Euler--Lagrange equation of a functional on periodic functions. A useful structural tool is a finite-depth version of the Plotnikov transform, with a zero-mode correction that is absent in infinite depth, originally introduced in \cite{plotnikov1992nonuniqueness} in the context of solitary waves and subsequently used in \cite{buffoni2000regularity,buffoni2000sub}. We introduce this transform explicitly in Section \ref{SecPlotnikov}.\\
We compare our results to the stability theorems proved in \cite{constantin2007stability}, where  the formal stability of a height-function variational problem, including the irrotational case, is studied. The constrained problems presented in this work differ from \cite{constantin2007stability} in that the wavelength and relative mass flux are fixed and the admissible perturbations preserve the average free-surface level, whereas in the principal result below we keep the conformal mean depth $h$ fixed and recover the mass flux from the Bernoulli constraint. Away from the trivial state, fixing conformal mean depth is not the same as fixing physical mean depth, see \eqref{relwd},
\begin{equation}
d-h=k[w\mathcal{C}_{kh}(w')]_{2\pi},
\end{equation}
which is quadratic in the wave amplitude. Consequently, the corresponding Hessians are restrictions to different constraint manifolds and their Morse indices need not agree. This explains why the shallow-water formal stability conclusion in \cite{constantin2007stability} does not contradict the fixed-$h$ instability established here.\\
The paper is structured as follows. In Section \ref{SecNotation}, we fix notation and distinguish variational stability from dynamical notions of stability. In Section \ref{SecPrelim}, we recall the variational formulation, its variational structure, and the local bifurcation theory. Section \ref{secstability} contains the second-variation calculation, the corrected stability analysis of the trivial solution, and the finite-depth Plotnikov transform. In Section \ref{SecPerturbation}, rather than truncating the transformed operator at first order, we compute the bifurcation curve to the order required for the Hessian sign and evaluate the second variation directly on the tangent to the non-trivial branch. Section \ref{SecConclusion} summarizes the results.\\
Our main result is summarized in the following theorem.
In the fixed-$h$ stability statement below, the second variation is taken with respect to $w$ alone, while $h$ and $Q$ are kept fixed. The mass flux $m$, which does not enter the $w$-equation or its $w$-Hessian, is then recovered from the Bernoulli constraint.
\begin{theorem}
    \label{mainthmintro}
Fix $h,k>0$ and consider variations in the mean-zero free-surface variable $w$ with $h$ fixed. At each critical point, $Q$ (equivalently $\mu$) and $k$ are held fixed when the $w$-Hessian is evaluated. The trivial solution is strictly variationally stable for
\begin{equation}
\mu>\mu_1^*=\frac{\tanh(kh)}{k},
\end{equation}
is variationally degenerate at $\mu=\mu_1^*$, and is variationally unstable for $\mu<\mu_1^*$. Moreover, the primary even non-trivial branch bifurcating from $(\mu_1^*,0)$ has a negative direction of the second variation for every sufficiently small non-zero amplitude. Hence the non-trivial waves are variationally unstable in the fixed-$h$ problem. On the unrestricted periodic space there is, in addition, the neutral translation mode generated by horizontal shifts of a non-trivial wave.
\end{theorem}


\section{Notation}\label{SecNotation}
Before we start with the main text, let us fix some notation. Let  $L^2_{2\pi}(\mathbb{R})$ denote the space of real-valued, $2\pi$-periodic square integrable functions with standard inner product
\begin{equation}
    \langle f,g\rangle =  \int_0^{2\pi} f(x) g(x)\, dx. 
\end{equation}
Any element $f\in L^2_{2\pi}(\mathbb{R})$ can be expanded in a Fourier series
\begin{equation}
f(x) = \sum_{n\in\mathbb{Z}} \hat{f}_n e^{\ri n x},    
\end{equation}
for the Fourier coefficients
\begin{equation}
    \hat{f}_n = \frac{1}{2\pi} \int_{0}^{2\pi} f(x) e^{-\ri n x}.
\end{equation}
We also write 
\begin{equation}
[f]_{2\pi}:=\frac{1}{2\pi}\int_{0}^{2\pi}f(x)\, dx,
\end{equation}
for the mean value of a function $f\in L^2_{2\pi}(\mathbb{R})$. The $s$-order Sobolev space is denoted as 
\begin{equation}
    H^{s}_{2\pi}(\mathbb{R}):=\left\{f\in L^2_{2\pi}(\mathbb{R}):\sum_{n\in\mathbb{Z}}(1+n^2)^s|\hat{f}_n| ^2 <\infty \right\},
\end{equation}
while the space of mean-free  $H^{s}$-functions is denoted as 
\begin{equation}\label{subscircle}
    H^{s}_{2\pi,\circ}(\mathbb{R}):=\left\{f\in H^{s}_{2\pi}(\mathbb{R}):[f]_{2\pi}=0\right\}.
\end{equation}
We denote the space of $2\pi$-periodic, $n$-times continuously  differentiable  functions with $\alpha$-H\"older continuous derivative as $C^{n,\alpha}_{2\pi}(\mathbb{R})$. Similar to the subscript $\circ$ in \eqref{subscircle} indicates functions with zero mean, the subscript $e$ indicates even functions.\\
Throughout, we write
\begin{equation}
    u' = \frac{du}{dx},
\end{equation}
for the $x$-derivative of a function to ease notation and sometimes abbreviate $d/dx$ as $\partial_x$. We denote the  nonzero integers by \begin{equation}
    \mathbb{Z}^* = \mathbb{Z}\setminus \{0\}.  
\end{equation}
We denote the strip of width $D$ as 
\begin{equation}
    \mathcal{S}_D=\{(x,y)\in\mathbb{R}^2 : -D<y<0\}.
\end{equation}

\begin{definition}\label{linstab}
Let $\mathcal{H}$ be a Hilbert space and let $\Lambda:\mathcal{H}\to\mathbb{R}$ be a $C^2$-functional. A critical point $w$ of $\Lambda$ is called \textit{strictly variationally stable} on a specified class of admissible variations $\mathcal{V}\subset\mathcal{H}$ if
\begin{equation}\label{deltadeltaJ}
    \delta^2\Lambda(w)(u,u)>0,
\end{equation}
for every nonzero $ u\in\mathcal{V}$. It is called \textit{weakly variationally stable} if the quadratic form is non-negative on $\mathcal{V}$, and \textit{variationally unstable} if there exists $u\in\mathcal{V}$ such that $\delta^2\Lambda(w)(u,u)<0$.
\end{definition}

\begin{remark}
We stress that the above notion of stability is related to a steady variational problem and should not be identified with linear stability of the time-dependent Cauchy problem. We further stress that the above notion of stability depends on the admissible class of functions. In particular, for a non-trivial periodic wave on the full periodic function space, horizontal translation invariance trivially produces a neutral direction. Thus strict positivity can only be expected after factoring out translations or after imposing a symmetry condition, for example by restricting to even perturbations. The instability result proved in the main text below is, however, stronger: it is produced by an even negative direction and therefore persists also on the full periodic space.
\end{remark}

\section{Preliminaries}\label{SecPrelim}

In this section, we recall the reformulation of the irrotational water wave problem in finite depth as a one-dimensional pseudo-differential equation. We summarize its variational structure as the Euler--Lagrange equation of a functional and state the local existence theory of small-amplitude solutions as proven in \cite{constantin2016global}. 

\subsection{Formulation of the Problem}

The governing equations for irrotational traveling water waves in finite depth can be reformulated as the following one-dimensional pseudo-differential equation,
\begin{equation}\label{equ}
\mu \mathcal{C}_{kh}(w')=\frac{w}{k}+w\mathcal{C}_{kh}(w')+\mathcal{C}_{kh}(ww')-[w\mathcal{C}_{kh}(w')]_{2\pi},
\end{equation}
for an unknown function $w\in H^1_{2\pi,\circ}(\mathbb{R})$ and an unknown scalar $h>0$, called the \textit{conformal mean depth}. The constant
\begin{equation}\label{defmu}
    \mu = \frac{Q-2gh}{g},
\end{equation}
contains constant hydraulic head $Q$ and the gravitational constant $g$, while $k>0$ is the wave number associated with the physical horizontal period. The pseudo-differential operator $\mathcal{C}_{D}: L^2_{2\pi,\circ}(\mathbb{R})\to L^2_{2\pi,\circ}(\mathbb{R})$, for any $D>0$, is called \textit{Hilbert transform for the strip of width $D$}, and acts on harmonic basis functions as 
a Fourier multiplier, \begin{equation}\label{defC}
    \mathcal{C}_{D}(e^{\ri n x}) = -\ri \coth(nD) e^{\ri n x},\quad n\neq 0. 
\end{equation}
A reformulation of the two-dimensional water wave problem with constant vorticity $\Upsilon$ as pseudo-differential equations has been derived in \cite{constantin2016global,constantin2011steady}, where equation \eqref{equ} is obtained in the irrotational case $\Upsilon = 0$. Indeed, writing
\begin{equation}\label{defv}
v = w +h,    
\end{equation}
any non-singular solution to \eqref{equ} that satisfies the constraint 
\begin{equation}\label{constr}
\left(\frac{m}{kh}\right)^2=\left[(Q-2gv)\Big((v')^2+\G{kh}(v)^2\Big)\right]_{2\pi},
\end{equation}
where $m$ is the mass flux,
is equivalent to the following system of equations and constraints
\begin{align}
& \left(\frac{m}{kh}\right)^2=(Q-2gv)\big((v')^2+\G{kh}(v)^2\big)\label{Bernw}\\
& [v]_{2\pi}=h\\
& v(x)>0 \text{ for } x\in[0,2\pi]\\
& \text{ the mapping } x\mapsto \left(\frac{x}{k}+\mathcal{C}_{kh}(v-h)(x),v(x)\right) \text{ is injective on }\mathbb{R}\\
& (v'(x))^2+\G{kh}(v)(x)^2\neq 0 \text{ for all } x\in[0,2\pi],
\end{align}
where the pseudo-differential operator $\mathcal{G}_{D}: L^2_{2\pi}(\mathbb{R})\to L^2_{2\pi}(\mathbb{R})$, called \textit{Dirichlet-to-Neumann operator}, acts on harmonic basis functions as 

\begin{equation}\label{defG}
    \mathcal{G}_{D}(e^{\ri n x}) =
      \begin{cases}
    n \coth(nD)e^{\ri n x},& n\neq 0,\\[0.1cm]
    D^{-1},& n=0,
    \end{cases}
\end{equation}
where the value at $n=0$ is understood by continuity. The operator is the Dirichlet-to-Neumann operator for the strip: it maps boundary Dirichlet data to the corresponding normal derivative of its harmonic extension \cite{constantin2016global}. The operators \eqref{defC} and \eqref{defG} are related through the identity
\begin{equation}\label{GCrelation}
    \mathcal{G}_D(u) = \frac{1}{D}[u]_{2\pi} + \mathcal{C}_{D}(u'). 
\end{equation}
Equation \eqref{Bernw} is a reformulation of Bernoulli's law for the traveling wave solution, i.e., a quasi-hydrostatic energy balance and we define the Bernoulli constant as
\begin{equation}\label{defB}
    B = \frac{1}{g}\left(\frac{m}{kh}\right)^2.
\end{equation}
Any solution to equation (\ref{equ}) that describes a regularly parametrized surface necessarily satisfies the inequality
\begin{equation}\label{inequ}
\frac{Q}{2g}> w+h
\end{equation}
as it can be seen from (\ref{Bernw}).\\
The advantage of the equation (\ref{equ}) as compared to (\ref{Bernw}) is that it has a variational formulation. Because the admissible $w$-variations are constrained to have zero mean, the Euler--Lagrange equation is understood after projection onto the mean-zero subspace. Namely, the first variation with respect to  $w\in H^1_{2\pi, \circ}(\mathbb{R})$, of the functional 
\begin{equation}\label{functional}
\begin{split}
\Lambda (w,h)=&\int_{0}^{2\pi} \bigg(Qv-gv^2\bigg)\left(\frac{1}{k}+\mathcal{C}_{kh}(w')\right)+\frac{m^2}{kh}\, dx
\end{split}
\end{equation}
satisfies $P_{\circ}(\delta\Lambda/\delta w)=2g\,\mathcal{F}(\mu,w)$, with $P_{\circ}f=f-[f]_{2\pi}$ and $\mathcal{F}$ defined below in \eqref{defF}, which gives the mean-zero Euler--Lagrange equation  \eqref{equ}. The first variation of \eqref{functional} with respect to $h$ recovers the constraint \eqref{constr}. We remark that the original variational problem in \cite{constantin2016global} was posed on the H\"{o}lder space $C^{2,\alpha}_{2\pi}$, which we relax here as $w\in H^1_{2\pi, \circ}(\mathbb{R})$ is sufficient to formulate the governing equations as critical points of \eqref{functional}. 
\begin{remark}
Since the mass flux $m$ does not appear in equation \eqref{equ}, the constraint \eqref{constr} defines the physical quantity $m$ in terms of the solution $w$  and the constants $Q$, $k$ and $h$. In the case of waves with non-vanishing vorticity as described in \cite{constantin2016global}, the mass flux, along with the constant vorticity $\Upsilon$, also appears in the first variation of $\Lambda $ and hence cannot be treated independently from $w$.  \end{remark}

\begin{remark}
For waves of infinite depth, an equation similar to (\ref{equ}) for a parametrization of the wave profile can be derived \cite{buffoni2003analytic,toland2000stokes,toland2002pseudo}:
\begin{equation}\label{equinf}
\tilde{\mu}\C{}(w')=w+w\C{}(w')+\C{}(ww'),
\end{equation}
where the wave number is normalized to $k=1$ and  $\mathcal{C}$ denotes the periodic Hilbert transform. On mean-zero periodic functions it is defined by the Fourier multiplier
\begin{equation}
    \mathcal{C}(e^{\ri n x}) = -\ri  \operatorname{sgn} (n) e^{\ri n x},\qquad n\neq0 .
\end{equation}
In fact,  (\ref{equinf}) can be viewed, at least formally, as the limiting equation of (\ref{equ}) as $h\to\infty$. A Lagrangian for equation (\ref{equinf}) is given by
\begin{equation}\label{funcinf}
\Lambda _{\infty}(w)=\int_{0}^{2\pi}\tilde{\mu} w\C{}(w')-w^2\big(1+\C{}(w')\big)\, dx.
\end{equation}
\end{remark}

\begin{remark}
Let us comment on how solutions to \eqref{equ} with the constraint \eqref{constr}, or, equivalently, equation \eqref{Bernw}, relate to the profile of physical water waves and how the conformal mean depth can be interpreted for irrotational waves. \\
The free surface of the irrotational traveling wave is recovered as 
\begin{equation}
 S = \left\{\left(\frac{x}{k}+\mathcal{C}_{kh}(v-h)(x),v(x)\right)\, ,\, x\in [0,2\pi] \right\}.
\end{equation}
 Since $\mathcal{C}_{kh}(v-h)$ is $2\pi$-periodic, the first component of this parametrization increases by exactly $2\pi/k$ when $x$ increases by $2\pi$. Hence the physical horizontal wavelength is $L=2\pi/k$, i.e. $k=2\pi/L$, see \cite{constantin2016global}. In the case of irrotational traveling waves, the conformal mean depth $h$ coincides with the mean depth as defined in \cite{amick1981periodic},
\begin{equation}\label{conmeand}
h=\frac{m}{c},
\end{equation}
where $c$ denotes the wave speed, or, equivalently, the constant speed of the traveling frame, with the sign convention for $m$ used in \cite{amick1981periodic}).\\
The free surface associated to a non-singular solution to (\ref{equ})-(\ref{constr}) is given by the graph of a function if and only if 
\begin{equation}\label{graph}
\frac{1}{k}+\mathcal{C}_{kh}(w')(x)> 0,
\end{equation}
for all $x\in [0,2\pi]$. If (\ref{graph}) holds true, denoting $I(x)=x/k+\mathcal{C}_{kh}(w)(x)$, we can write the free surface in terms of the parametrization $v$ as a graph $$\eta(x):=v(I^{-1}(x)).$$
Note that $I$ is a smooth bijection from $[0,2\pi]$ onto $[0,L]$. This allows us to relate the mean depth $d$, defined such that $\eta-d$ is mean-free over a period $L$, to the  conformal mean depth $h$ as
\begin{equation*}
\begin{split}
d &=\frac{1}{L}\int_0^L\eta(x)\, dx\\
& =\frac{1}{L}\int_0^Lv(I^{-1}(x))\, dx\\
&=\frac{1}{L}\int_{0}^{2\pi}v(\xi)\left(\frac{1}{k}+\mathcal{C}_{kh}(w')(\xi)\right), d\xi\\
&=\frac{2\pi h}{kL}+\frac{2\pi}{L}[w\mathcal{C}_{kh}(w')]_{2\pi}\\[0.3cm]
&=h+k[w\mathcal{C}_{kh}(w')]_{2\pi},
\end{split}
\end{equation*}
where we have transformed the integral according to $ x=I(\xi)$ and used the fact that $v=w+h$ with $[w]_{2\pi}=0$. This means that the average appearing in equation \eqref{equ} can be interpreted as
\begin{equation}\label{relwd}
[w\mathcal{C}_{kh}(w')]_{2\pi}=\frac{d-h}{k},
\end{equation}
i.e., the deviation of the physical mean depth from the conformal mean depth.
\end{remark}

\subsection{Existence of Solutions and Global Bifurcation Analysis}

Let us recall the existence theory of solutions to equation \eqref{equ}. Clearly, $w=0$ defines a solution to \eqref{equ} for any value of $\mu$, $k$ and $h$. By the constraint \eqref{constr}, the hydraulic head $Q$, the conformal mean depth $h$ and the mass flux $m$ at the trivial solution are related as
\begin{equation}\label{Qtriv}
Q=2gh+\frac{m^2}{h^2}.
\end{equation}
Indeed, at $v=h$ one has $\mathcal{G}_{kh}(h)=1/k$, so that \eqref{constr} gives $(m/(kh))^2=(Q-2gh)/k^2$. In view of \eqref{conmeand}, $m=ch$ at the trivial state, hence $Q-2gh=c^2$ and $\mu=c^2/g$. The bifurcation value \eqref{munstar} for $n=1$ therefore recovers $c^2=(g/k)\tanh(kh)$.
We denote the surface of trivial solutions, parametrized by the parameter $\mu$ or, equivalently, the hydraulic head, and the conformal mean depth as \begin{equation}
    \mathcal{K}_{\rm triv} = \{(h,\mu,0),h\in\mathbb{R}^{+},\mu\in\mathbb{R}\} \subset \mathbb{R}^+\times\mathbb{R}\times L^2_{2\pi,\circ}(\mathbb{R}). 
\end{equation}

\begin{theorem}\label{thmexistence}
Let $h,k>0$ be given and denote 
\begin{equation}\label{munstar}
\mu_{n}^{*} = \frac{\tanh(nkh)}{nk},\quad n\in\mathbb{N}. 
\end{equation}
For any $\mu\in\mathbb{R}\setminus\{\mu_{n}^*, n\in\mathbb{N}\}$, there exists a neighborhood in $\mathbb{R}^+\times \mathbb{R}\times C^{2,\alpha}_{2\pi,e}(\mathbb{R})$ of the point $(h,\mu,0)$ on $\mathcal{K}_{\rm triv}$ in which the only solutions of \eqref{equ} are those on $\mathcal{K}_{\rm triv}$.
For each integer $n\geq 1$,  there exists a continuous curve 
\begin{equation}
\mathcal{K}_{n} = \{(\mu(\varepsilon),w_\varepsilon):\varepsilon\in\mathbb{R}\},
\end{equation}
of solutions $w_\varepsilon$ to \eqref{equ} in the space $\mathbb{R}\times C^{2,\alpha}_{2\pi,e}(\mathbb{R})$, such that the following properties hold:
\begin{enumerate}
    \item The curve of non-trivial solutions bifurcates from the trivial branch as $(\mu(0),w_0)=(\mu_n^*,0)$, where $\mu_n^*$ is given by \eqref{munstar}.
    \item The non-trivial solution admits the leading-order Taylor expansion
    \begin{equation}\label{veps}
   w_\varepsilon(x) = \varepsilon\cos(nx) +o(\varepsilon),\quad 0<|\varepsilon|<\varepsilon_0
\end{equation}
for $\varepsilon_0$ sufficiently small. 
\item There exists a neighborhood $\mathcal{W}_{n}$ of $(\mu_n^*,0)$ in $\mathbb{R}\times C^{2,\alpha}_{2\pi,e}(\mathbb{R})$ and $\varepsilon$ sufficiently small such that
\begin{equation}
    \{ (\mu,w) \in \mathcal{W}_{n}: w\neq0 \text{ and } \eqref{equ} \text{ holds }  \} = \{(\mu(\varepsilon), w_\varepsilon): 0<|\varepsilon|<\varepsilon_0\}
\end{equation}
\item The corresponding surface variable $v_\varepsilon=h+w_\varepsilon$ satisfies $\mu(\varepsilon)/2>w_\varepsilon$.
\item The curve $\mathcal{K}_n$ has a real-analytic re-parametrization locally around each of its points. 

\end{enumerate}

\end{theorem}
A proof of the above theorem, based on the Crandall--Rabinowitz bifurcation theorem, is given in \cite{constantin2016global}. Theorem \ref{thmexistence} will serve as the basis of the stability analysis of the non-trivial solution branch for small amplitudes in the following sections.  The definition of $\mu_n^*$ in \eqref{munstar} is the classical dispersion relation for traveling water waves under the influence of gravity \cite{constantin2011nonlinear}.

\section{Stability Analysis of Traveling Waves}\label{secstability}

In this section, we will investigate stability properties of critical points of the functional (\ref{functional}). We first analyze the stability of the trivial solution  and compare it to the trivial solution in infinite depth in Subsection \ref{SecStabilityTrivial}. Then, in Subsection \ref{SecPlotnikov}, we employ a finite-depth analog of the Plotnikov transform to rewrite the second variation of $\Lambda$  as a diagonal quadratic form plus  a multiplication term and an explicit finite-rank zero-mode correction.\\
The second variation of (\ref{functional}) with respect to $w$ for $u\in H^{1}_{2\pi,\circ}(\mathbb{R})$ is given by
\begin{small}
\begin{equation}\label{secderww}
\begin{split}
&\frac{\delta^2\Lambda }{\delta w^2}  (w,h)(u,u)  = \left. \frac{d^2}{d\varepsilon^2}\right|_{\varepsilon = 0} \Lambda (w+\varepsilon u,h)\\
  & = \left. \frac{d^2}{d\varepsilon^2}\right|_{\varepsilon = 0} \int_{0}^{2\pi} \bigg(Q(w+\varepsilon  u+h)-g(w+\varepsilon u+h)^2\bigg)\left(\frac{1}{k}+\mathcal{C}_{kh}(w') + \varepsilon \mathcal{C}_{kh}(u') \right)+\frac{m^2}{kh}\, dx \\
  & = 2\int_{0}^{2\pi} (Q-2gv)u\mathcal{C}_{kh}(u')-g\left(\frac{1}{k}+\mathcal{C}_{kh}(w')\right)u^2 \, dx. 
\end{split}
\end{equation}
\end{small}

\begin{remark}\label{remconstraints}
In the principal stability analysis we vary $w$ and keep the conformal mean depth $h$ fixed for the evaluation of the Hessian in the $w$-direction. More precisely, the admissible tangent space is $H^1_{2\pi,\circ}(\mathbb{R})$, and $h$, $Q$ (equivalently $\mu$), and $k$ are held fixed during this second variation. The parameter $m$ does not enter the $w$-equation or the $w$-Hessian and is recovered, for a physical wave, from \eqref{constr}. This choice of constraint is important. By \eqref{relwd}, physical mean depth and conformal mean depth differ on a non-trivial wave by the quadratic term $k[w\mathcal{C}_{kh}(w')]_{2\pi}$; consequently, fixing $h$ is not equivalent to fixing physical mean depth along the non-trivial branch. Likewise, the height-function problem in \cite[Section 4]{constantin2007stability}, which fixes relative mass flux and mean surface level, is a different constrained second-variation problem. We discuss variations in $h$ separately for the trivial solution below.
\end{remark}

\begin{theorem}
    \label{mainthm}
Let $h,k>0$ be fixed and let $\mu_1^*=\tanh(kh)/k$. The $w$-Hessian below is evaluated with $h$, $Q$ (equivalently $\mu$), and $k$ fixed. For mean-zero variations in $w$, the trivial solution $w=0$ is strictly variationally stable for
\begin{equation}\label{condmu}
    \mu> \frac{\tanh(kh)}{k},
\end{equation}
weakly variationally stable and degenerate at $\mu=\mu_1^*$, and variationally unstable for $\mu<\mu_1^*$. At $\mu=\mu_1^*$ the kernel is spanned by $\cos x$ and $\sin x$.\\
Let $(\mu(\varepsilon),w_\varepsilon)$ be the primary even branch with
$w_\varepsilon=\varepsilon\cos x+o(\varepsilon)$. Then, for every finite $h,k>0$ and all sufficiently small $\varepsilon\neq0$, there exists an even mean-zero variation $u_\varepsilon$ such that
\begin{equation}
\frac{\delta^2\Lambda}{\delta w^2}(w_\varepsilon,h)
(u_\varepsilon,u_\varepsilon)<0.
\end{equation}
Thus the non-trivial primary branch is variationally unstable in the fixed-$h$ problem.
\end{theorem}

\begin{figure}
    \centering
    \includegraphics[width=0.8\linewidth]{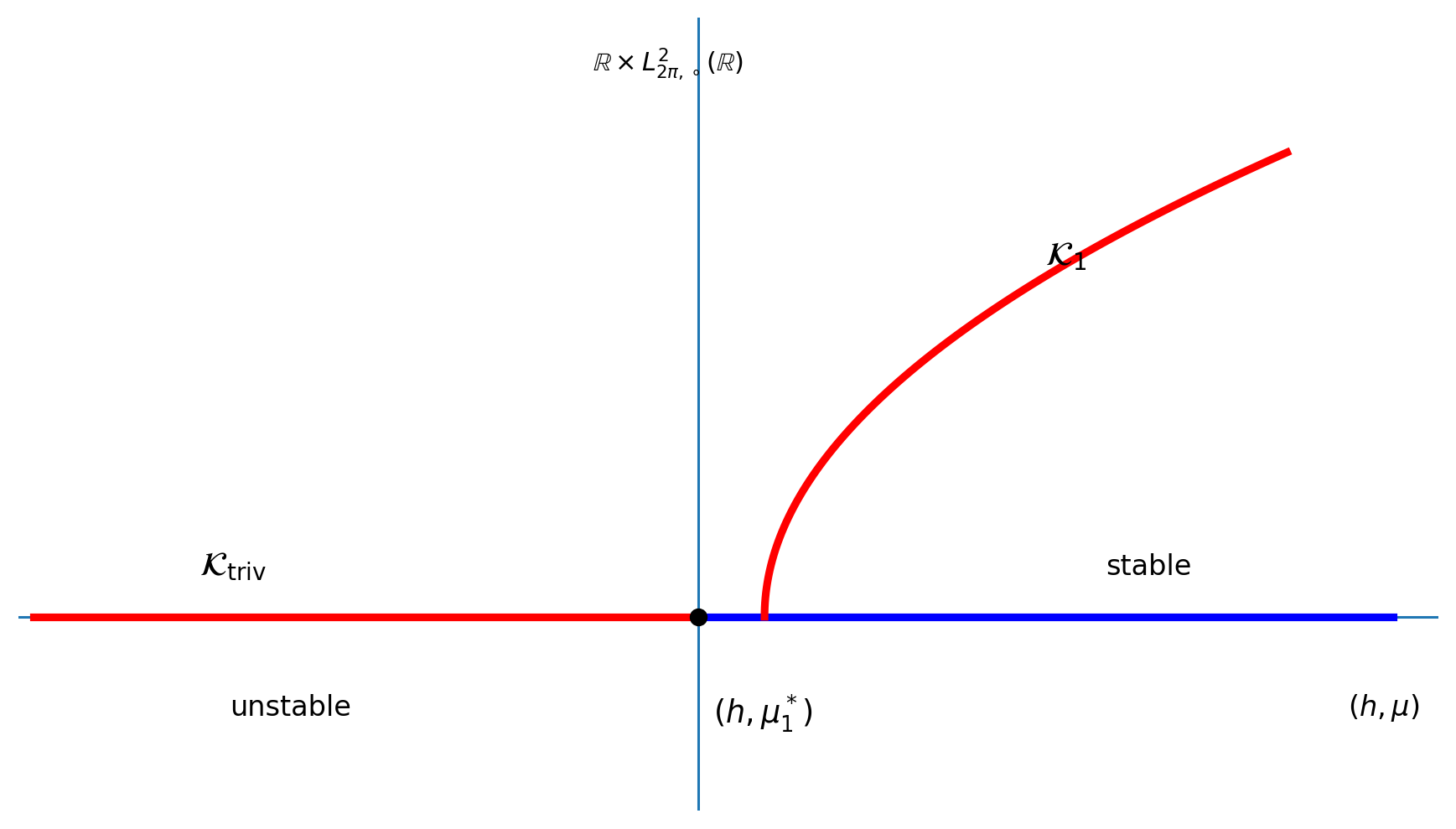}
    \caption{Schematic bifurcation diagram of solutions to \eqref{equ}, with $\mu$ increasing to the right: the trivial branch $\mathcal{K}_{\rm triv}$ is unstable for $\mu<\mu_1^*$ (red), degenerate at $\mu=\mu_1^*$, and strictly variationally stable for $\mu>\mu_1^*$ (blue). Since $\mu_2>0$, the primary non-trivial branch $\mathcal{K}_1$ bends toward larger $\mu$ and is variationally unstable for sufficiently small non-zero amplitudes (red curve).}
    \label{plot_bifurcation}
\end{figure}

The first assertion follows immediately from the Fourier representation of the Hessian and is proved in Subsection \ref{SecStabilityTrivial}. The instability of the primary non-trivial branch is proved in Section \ref{SecPerturbation}. Figure \ref{plot_bifurcation} shows a schematic depiction of the stability of solutions to \eqref{equ} around the first bifurcation point. The remainder of this section  is  devoted to the preparation of  the proof of Theorem \ref{mainthm}.

\subsection{Stability of the Trivial Solution}
\label{SecStabilityTrivial}

Before we turn to the investigation of stability for general solutions, let us examine $\delta^2\Lambda/\delta w^2$ at the trivial solution $w=0$. This proves the first statement of Theorem \ref{mainthm}. The second variation simplifies to
\begin{equation}\label{secondvartriv}
 \begin{split}
\frac{\delta^2\Lambda }{\delta w^2}(0,h)(u,u) & = 2\int_{0}^{2\pi}(Q-2gh)u\mathcal{C}_{kh}(u')-\frac{g}{k}u^2\, dx\\[0.3cm]
&=2\int_{0}^{2\pi}g\mu u\mathcal{C}_{kh}(u')-\frac{g}{k}u^2\, dx\\[0.3cm] 
& = 4\pi g\sum_{n\in\mathbb{Z}^*} \left(\mu n \coth(nkh)-\frac{1}{k}\right)|\hat{u}_n|^2.
\end{split}
\end{equation}
The function $n\mapsto n\coth(nkh)$ is increasing for $n\geq1$. Hence the smallest Fourier coefficient in \eqref{secondvartriv} is the $|n|=1$ coefficient. It follows that the quadratic form is strictly positive if and only if
\begin{equation}
\mu>\frac{\tanh(kh)}{k}.
\end{equation}
At $\mu=\mu_1^*$ it is non-negative with kernel $\operatorname{span}\{\cos x,\sin x\}$, while for $\mu<\mu_1^*$ the first harmonic gives a negative direction. In particular, the bifurcation point is degenerate rather than strictly unstable.\\
At a higher bifurcation point $\mu_n^*$, $n>1$, we have
\begin{equation}\label{secondvarzero}
 \begin{split}
   \frac{\delta^2\Lambda }{\delta w^2}(0,h)(u,u) =  \frac{4\pi g}{k}\sum_{j\in\mathbb{Z}^*} \left( \frac{j \coth(jkh)}{n\coth(nkh)} -1\right)|\hat{u}_j|^2.
\end{split}
\end{equation}
Thus the negative subspace is spanned by
\begin{equation}
\{\cos(jx),\sin(jx):1\leq j<n\},
\end{equation}
and has real dimension $2(n-1)$, while $\cos(nx)$ and $\sin(nx)$ span the two-dimensional kernel.\\
As mentioned in Remark \ref{remconstraints}, the main theorem concerns the fixed-$h$ Hessian. For comparison only, we now inspect the full variation in $(w,h)$ at the trivial solution. In this auxiliary calculation the parameters $Q$, $m$, and $k$ in the functional \eqref{functional} are held fixed while $(w,h)$ is varied. Provided that the fixed-$h$ condition \eqref{condmu} holds, the trivial flow is strictly variationally stable also with respect to variations in the conformal mean depth if, in addition,
 \begin{equation}\label{condmgh}
m^2 > gh^3.
\end{equation}
The second partial variation with respect to $h$ at the trivial solution is
\begin{equation}
\begin{split}
\frac{\delta^2\Lambda}{\delta h^2}(0,h) & = \frac{\partial^2}{\partial h^2} \int_0^{2\pi} \frac{Qh-gh^2}{k} +\frac{m^2}{kh} \, dx \\
& = \int_0^{2\pi} -\frac{2g}{k} + \frac{2m^2}{kh^3}\, dx.
\end{split}
\end{equation}
Since $\big(\delta^2\Lambda /\delta w\delta h\big)(0,h)(u,\nu)=0$ for all $u\in H_{2\pi,\circ}^1(\mathbb{R})$ and $\nu\in\mathbb{R}$, the full second variation reduces to
\begin{equation}
\begin{split}
\delta^2\Lambda (0,h)((u,\nu),(u,\nu)) =2\int_{0}^{2\pi}(Q-2gh)u\mathcal{C}_{kh}(u')-\frac{g}{k}u^2
+\left(\frac{m^2}{kh^3}-\frac{g}{k}\right)\nu^2\, dx.
\end{split}
\end{equation}
At the trivial solution $d=h$, and \eqref{conmeand} gives $m=ch$. Therefore \eqref{condmgh} is equivalent to the strict inequality
\begin{equation}\label{condtrivsol}
c^2>gd.
\end{equation}
The equality $c=\sqrt{gd}$ is the classical long-wave critical speed. By contrast, the linear phase speed at a finite-wavenumber periodic bifurcation satisfies
\begin{equation}
c_{\rm per}^2=\frac{g}{k}\tanh(kd)<gd,
\end{equation}
with equality recovered only in the long-wave limit $k\downarrow0$. Small-amplitude solitary waves bifurcate from the critical long-wave speed $c=\sqrt{gd}$ and travel on the supercritical side. This places \eqref{condtrivsol} in the broader steady-water-wave context.\\
We stress again that allowing $h$-variations changes the variational problem. The condition \eqref{condmgh} is an additional requirement for positivity of the full $(w,h)$ Hessian and should not be conflated with the fixed-$h$ stability condition \eqref{condmu}.

\begin{figure}[h!]
    \centering
\includegraphics[width=0.7\linewidth]{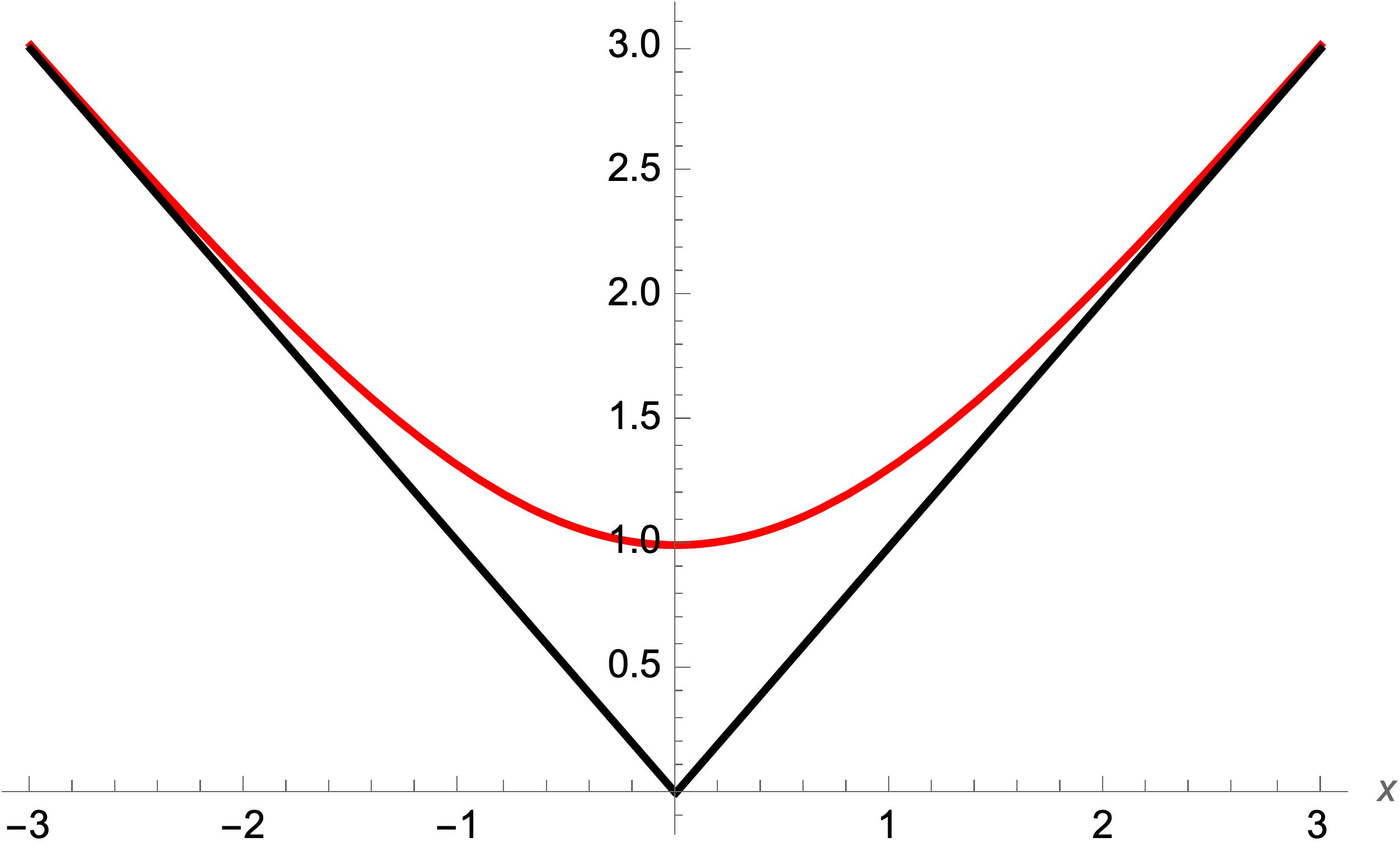}
    \caption{The symbols $\xi\coth(D\xi)$ (red, $D=1$) and $|\xi|$ (black). The finite-depth symbol extends continuously to the value $D^{-1}$ at $\xi=0$, which is the zero-mode eigenvalue of the full Dirichlet-to-Neumann operator $\mathcal{G}_D$. On the mean-zero $2\pi$-periodic space, however, the spectrum of $\mathcal{C}_D\partial_x$ is sampled only at nonzero integers and its smallest eigenvalue is $\coth D$, while the corresponding infinite-depth smallest eigenvalue is $1$.}
    \label{figsymbol}
\end{figure}

\begin{remark}
We contrast the stability of the trivial solution for finite-depth traveling waves to the trivial solution in infinite depth. The second derivative of $\Lambda _{\infty}$ as defined in \eqref{funcinf} at the trivial solution $w=0$ is given by
\begin{equation}
\delta^2\Lambda _{\infty}(0)(u,u)=2\int_{0}^{2\pi}\left(\tilde{\mu}\C{}(u')-u\right)u\, dx.
\end{equation}
On the full periodic space $H^1_{2\pi}(\mathbb{R})$, this quadratic form cannot be positive definite for any choice $\tilde{\mu}\geq0$, because the constant Fourier mode is an eigenfunction of $\tilde{\mu}\mathcal{C}\partial-1$ with eigenvalue $-1$. This is the setting in which the infinite-depth critical points are saddles, as pointed out in \cite{buffoni2003surface}. This observation should not be confused with the fixed-$h$ finite-depth problem above, where the admissible $w$-variations are mean-zero and the constant mode is excluded. On that mean-zero periodic space, the finite-depth multiplier $\mathcal{C}_D\partial_x$ has eigenvalues $|n|\coth(|n|D)$, $n\neq0$, with smallest value $\coth D$, whereas the corresponding infinite-depth multiplier has eigenvalues $|n|$, with smallest value $1$. The continuous symbol $\xi\coth(D\xi)$ tends to $D^{-1}$ as $\xi\to0$; this value is the zero-mode eigenvalue of the full Dirichlet-to-Neumann operator $\mathcal{G}_D$, not an eigenvalue of $\mathcal{C}_D\partial_x$ on the mean-zero space. Figure \ref{figsymbol} illustrates this distinction.\\
\end{remark}

\subsection{Rewriting the second variation $\delta^2\Lambda /\delta w^2$ by the Plotnikov transform}\label{SecPlotnikov}

In this subsection, we pull back the second variation $\delta^2\Lambda /\delta w^2$ by a linear  transformation adapted from 
\cite{plotnikov1992nonuniqueness}. In finite depth one has to retain a zero Fourier mode generated by products of mean-zero functions. This produces an explicit finite-rank correction to the infinite-depth Plotnikov formula. Since all variations considered below have zero mean, we formulate  the holomorphic extension  directly on the mean-zero space. Throughout this subsection $w$ is a regular solution, for instance $w\in C^{2,\alpha}_{2\pi}$ as supplied by the local bifurcation theory. Consequently multiplication by $w'$ and by $1/k+\mathcal{C}_D(w')$ is bounded on $H^1$, so the map $\mathcal{P}$ defined below is a bounded linear map $H^1_{2\pi,\circ}\to H^1_{2\pi,\circ}$. The identity in Proposition \ref{secderwwPlot} is a pullback formula on the range of $\mathcal P$ and no global invertibility of $\mathcal P$ is required for the instability proof in Section \ref{SecPerturbation}.\\
Henceforth, we write  $z=x+\ri y$. For $u\in H^1_{2\pi,\circ}(\mathbb{R})$ with Fourier expansion
\begin{equation}
    u(x)=\sum_{n\in\mathbb{Z}^*}\hat u _n e^{\ri n x},
\end{equation}
  {define
} \begin{equation}\label{holex}
\mathbf{E} _D[u] (z):= -\sum_{n\in\mathbb{Z}^*}\frac{e^{-nD}}{\sinh(nD)} \hat{u}_n e ^{\ri n z}.
\end{equation}
 Then $\mathbf{E}_D[u]$ is $2\pi$-periodic  and holomorphic in  $\mathcal{S}_D$, and
\begin{equation}\label{RealEzero}
    \Re\mathbf{E}_D[u](x-\ri D)=0,
\end{equation}
 while at the upper boundary
 \begin{equation}\label{relE}
    \mathbf{E}_D[u](x)=u(x)-\ri\mathcal{C}_D(u)(x).
\end{equation}
 Its real part is the harmonic  solution of the  Dirichlet problem with boundary values $0$ at $y=-D$ and $u$ at $y=0$. Explicitly,
\begin{equation}\label{solDir}
U(x,y)=\sum _{n\in\mathbb{Z}^*} \frac{\sinh\big(n(y+D)\big)}{\sinh(nD)} \hat u_ne ^{\ri n x}=\Re\mathbf{E}_D[u](x+\ri y).
 \end{equation}
 The restriction to mean-zero data is essential for the periodic holomorphic extension: for nonzero mean, the harmonic conjugate contains a term linear in $x$.
 To define the Plotnikov transform, write
\begin{equation}\label{defW}
W(z)=\frac{1}{k}+\ri\mathbf{E} _D[w'] (z).
\end{equation}
 For $x\in[0,2\pi]$,
 \begin{equation}\label{defW0}
W(x)=\frac{1}{k}+\mathcal{C}  {_D} (w')(x)+\ri w'(x)  .
 \end{equation}
For the regular graph solutions considered here, $W$ has no zeros in the closed strip. Indeed, $\Re W>0$ on the upper boundary by \eqref{graph}. On the lower boundary $\Im W=0$, and the Cauchy--Riemann equations therefore give $\partial_y\Re W=0$. If $\Re W$ attained a non-positive minimum in the interior, the strong maximum principle would force $\Re W$ to be constant, contradicting its strictly positive upper-boundary values. If such a minimum were attained on the lower boundary, Hopf's boundary point lemma would force a nonzero inward normal derivative, contradicting $\partial_y\Re W=0$ there. Hence $\Re W>0$ throughout the closed strip, and in particular $W$ has no zeros.

 \begin{definition}\label{defP}
  {For $u\in H^1_{2\pi,\circ}(\mathbb{R})$ and } fixed $W$ of the form \eqref{defW},  {define } the \textit{Plotnikov transform}
 \begin{equation}  \label{PlotnikovDef}
{\mathcal{P}[u]:=\Re\{W\mathbf{E}_D[u]\}
=u}\left({\frac{1}{k}+\mathcal{C}_D(w')}\right){+w'\mathcal{C}_D(u).
} \end{equation}
 \end{definition}
  {By skew-symmetry of $\mathcal{C}_D$, one has $[\mathcal{P}[u]]_{2\pi}=0$.
}

{The finite-depth analogue of the usual Hilbert-transform product identity contains a zero-mode correction.
} \begin{lemma} \label{lemFiniteDepthTricomi}
 Let   {$u,v\in H^1_{2\pi,\circ}(\mathbb{R})$. Then
}\begin{equation}{\label{idCD}
\mathcal{C}_D\big(u\mathcal{C}_D(v)+v\mathcal{C}_D(u)\big)
=\mathcal{C}_D(u)\mathcal{C}_D(v)-uv
-\big[\mathcal{C}_D(u)\mathcal{C}_D(v)-uv\big]_{2\pi}.
}\end{equation} 
 \end{lemma}
\begin{proof}
 {It is enough to prove the identity first for trigonometric polynomials and then pass to $H^1$ by density. } Expanding $u$ and $v$   {in Fourier series, the $n$th Fourier coefficient for $n\neq0$ is computed using
}\begin{equation}
{\coth(x+y)=\frac{1+\coth(x)\coth(y)}{\coth(x)+\coth(y)}.
}\end{equation} 
  {The terms for which one of the input frequencies is zero vanish because $u$ and $v$ have zero mean. Thus, for every $n\neq0$, the two sides of \eqref{idCD} have the same Fourier coefficient. The left-hand side has zero mean because $\mathcal{C}_D$ maps into the mean-zero space. Subtracting the mean on the right gives equality also in the zero Fourier mode} .
\end{proof}

  {For later use, set
} \begin{equation}  \label{defMD}
{M_D[u;w]:=}\big[{\mathcal{C}_D(u)\mathcal{C}_D(w')-uw'}\big]{_{2\pi}.
} \end{equation}
  {Applying Lemma \ref{lemFiniteDepthTricomi} with $v=w'$ gives the corrected multiplicative relation
}\begin{equation}{\label{propP}
\mathbf{E}_D[\mathcal{P}[u]]=W\mathbf{E}_D[u]+\ri M_D[u;w].
}\end{equation} 
  {The additional term is a purely imaginary constant. It disappears after differentiation, so that
} \begin{equation}\label{derRP}
\mathbf{E}_D[\mathcal{P}[u]']
= (W\mathbf{E}_D[u])'
=W\mathbf{E}_D[u']+W'\mathbf{E}_D[u] {.
} \end{equation}
Notice that the constant in \eqref{propP} is generally nonzero even though both $u$ and $w'$ have zero mean.
We now set $D=kh$. Bernoulli's law \eqref{Bernw} and \eqref{defB} give
\begin{equation}\label{BernW}
(Q-2gv)|W|^2=gB,
\end{equation}
which, after differentiation becomes \begin{equation} \label{derBernw}
{w} '|W|^2=  {B\Re}\left\{{\frac{W'}{W}}\right\}{.
} \end{equation}
  {Define the Plotnikov potential
}\begin{equation}{\label{defPhi}
\Phi:=\Im\left\{\frac{W'}{W}\right\}
+\frac{|W|^2}{B}\left(\frac{1}{k}+\mathcal{C}_{kh}(w')\right).
}\end{equation} 
\begin{proposition}\label{secderwwPlot}
For $u\in H^1_{2\pi,\circ}(\mathbb{R})$,  the second variation   {satisfies
}\begin{equation}{\label{CorrectedPlotnikov}
\begin{split}
\frac{\delta^2\Lambda}{\delta w^2}\big(\mathcal{P}[u],\mathcal{P}[u]\big)
={}&2gB\int_0^{2\pi}u\big(\mathcal{C}_{kh}(u')-\Phi u\big)\,dx\\
&+2gB\,M_{kh}[u;w]\int_0^{2\pi}\frac{w'}{|W|^2}
\big(\mathcal{C}_{kh}(u')-\Phi u\big)\,dx.
\end{split}
}\end{equation} 
  {Thus the usual diagonal-plus-potential Plotnikov form is supplemented in finite depth by an explicit finite-rank term. In the infinite-depth limit this correction disappears because the Hilbert-transform multiplier has modulus one on every nonzero Fourier mode} .
\end{proposition}

To prove Proposition \ref{secderwwPlot}, we   use the following elementary contour identity.
 \begin{lemma}\label{Imparts}
Let $F,G:\mathcal{S}_D\to\mathbb{C}$ be $2\pi$-periodic holomorphic functions. Assume   that
\begin{equation}\label{assRe}
\Re  F (x-\ri D)=\Re  G (x-\ri D)=0.
\end{equation}
Then
  \begin{equation}
{\int_0^{2\pi}\Im\{F^*G\}(x)\,dx
=2\int_0^{2\pi}\Re F(x)\Im G(x)\,dx.
}\end{equation} 
\end{lemma}
 \begin{proof}
  {Apply } Cauchy's   {integral theorem to } $FG$   {on the rectangle $(0,0)\to(2\pi,0)\to(2\pi,-D)\to(0,-D)$. The vertical contributions cancel by periodicity and the real parts vanish on } the lower edge .   {Taking imaginary parts gives
}\begin{equation}
{\int_0^{2\pi}\big(\Re F\,\Im G+\Re G\,\Im F\big)\,dx=0.
}\end{equation} 
 Since $\Im(F^*G)=\Re F\,\Im G-\Im F\,\Re G$, the assertion follows.
\end{proof}

\begin{proof}[Proof of Proposition \ref{secderwwPlot}]
For a mean-zero $q\in H^1_{2\pi,\circ}(\mathbb{R})$, integration by parts in \eqref{secderww} gives
 \begin{equation}\label{trans1}
  \begin{split}
\frac{\delta^2\Lambda}{\delta w^2}(q,q)
={}&\int_0^{2\pi}(Q-2gv)
\big(q\mathcal{C}_{kh}(q')-q'\mathcal{C}_{kh}(q)\big)\,dx\\
&+2g\int_0^{2\pi}q\left\{w'\mathcal{C}_{kh}(q)
-q\left(\frac{1}k+\mathcal{C}_{kh}(w')\right)\right\}\,dx.
\end{split} 
\end{equation}
 To make this step explicit, set
\begin{equation}
F=\frac{\mathbf{E}_{kh}[q]}{W},\qquad
G=\frac{\mathbf{E}_{kh}[q']}{W}.
\end{equation}
By \eqref{relE},
\begin{equation}
\Im\!\left(\mathbf{E}_{kh}[q]^*\mathbf{E}_{kh}[q']\right)
=-\big(q\mathcal{C}_{kh}(q')-q'\mathcal{C}_{kh}(q)\big).
\end{equation}
Using \eqref{BernW}, the first integral in \eqref{trans1} is therefore
\begin{equation}
-gB\int_0^{2\pi}\Im(F^*G)\,dx
=-2gB\int_0^{2\pi}\Re F\,\Im G\,dx,
\end{equation}
where Lemma \ref{Imparts} was used in the last equality. Writing $W=a+\ri b$, with $a=1/k+\mathcal{C}_{kh}(w')$ and $b=w'$, one also has
\begin{equation}
\Re F=\frac{aq-b\mathcal{C}_{kh}(q)}{|W|^2},
\end{equation}
so the second integral in \eqref{trans1} equals $-2g\int_0^{2\pi}|W|^2q\,\Re F\,dx$. Combining the two contributions gives
\begin{equation}  \label{trans2}
\frac{\delta^2\Lambda}{\delta w^2}(q,q)
=-2g\int_0^{2\pi}\Re\left\{\frac{\mathbf{E}_{kh}[q]}{W}\right\}
\left(B\Im\left\{\frac{\mathbf{E}_{kh}[q']}{W}\right\}+|W|^2q\right)\,dx.
\end{equation}
Now   {set $q=\mathcal{P}[u]$. From \eqref{propP},
} \begin{equation}  \label{firstfactorcorrected}
{\Re}\left\{{\frac{\mathbf{E}_{kh}[\mathcal{P}[u]]}{W}}\right\}
{=u+M_{kh}[u;w]\Re}\left\{{\frac{\ri}{W}}\right\}
{=u+M_{kh}[u;w]\frac{w'}{|W|^2}.
} \end{equation}
  {On the other hand, the constant correction in \eqref{propP} disappears after differentiation. Hence \eqref{derRP}, \eqref{derBernw}, and \eqref{defPhi} give
}\begin{equation}{\label{secondfactorcorrected}
B\Im\left\{\frac{\mathbf{E}_{kh}[\mathcal{P}[u]']}{W}\right\}
+|W|^2\mathcal{P}[u]
=-B\mathcal{C}_{kh}(u')+B\Phi u.
}\end{equation} 
  {Substituting \eqref{firstfactorcorrected} and \eqref{secondfactorcorrected} into \eqref{trans2} yields \eqref{CorrectedPlotnikov}} .
\end{proof}

 \section{Instability of Solutions along the Primary Bifurcation Branch}\label{SecPerturbation}

We now prove the second assertion of Theorem \ref{mainthm}. The argument is carried out directly in the variational formulation   {and therefore does not rely on the transformed expression in Proposition \ref{secderwwPlot}. In addition to the finite-rank zero-mode correction displayed there, a } first-order expansion of the Plotnikov potential   {alone would not determine a second-order eigenvalue correction when the first-order correction vanishes} . Indeed, an expansion of the form
\begin{equation}
\Phi_\varepsilon=\Phi_0+\varepsilon\Phi_1+o(\varepsilon)
\end{equation}
only yields
\begin{equation}
\mathcal{L}_\varepsilon=\mathcal{L}_0+\varepsilon\mathcal{L}_1+o(\varepsilon).
\end{equation}
If the first-order eigenvalue correction vanishes, then the $O(\varepsilon^2)$ part of the operator contributes at the same order as the usual term
$-\mathcal{L}_1\mathcal{L}_0^{-1}\mathcal{L}_1$ and cannot be omitted. \\
For fixed $h,k>0$, define
\begin{equation}\label{defF}
\begin{split}
\mathcal{F}(\mu,w):={}&\mu\mathcal{C}_{kh}(w')-\frac{w}{k}
-w\mathcal{C}_{kh}(w')-\mathcal{C}_{kh}(ww')
+[w\mathcal{C}_{kh}(w')]_{2\pi}.
\end{split}
\end{equation}
Then \eqref{equ} is equivalent to $\mathcal{F}(\mu,w)=0$. To ease notation in the following, we set 
\begin{equation}
T:=\mathcal{C}_{kh}\partial_x.
\end{equation}
On the mean-zero space $T$ is self-adjoint and
\begin{equation}
T(\cos(nx))=n\coth(nkh)\cos(nx),\qquad
T(\sin(nx))=n\coth(nkh)\sin(nx).
\end{equation}
Differentiating \eqref{functional} and projection onto mean-zero functions gives
\begin{equation}\label{gradF}
P_{\circ}\frac{\delta\Lambda}{\delta w}(w,h)=2g\,\mathcal{F}(\mu,w),
\end{equation}
where $P_{\circ}f=f-[f]_{2\pi}$. Consequently, at a critical point and for every mean-zero variation $u$,
\begin{equation}\label{HessF}
\frac{\delta^2\Lambda}{\delta w^2}(w,h)(u,u)
=2g\langle u,\mathcal{F}_w(\mu,w)u\rangle.
\end{equation}
Formula \eqref{HessF} is equivalent to \eqref{secderww} but is particularly convenient along a bifurcating solution curve.

\subsection{Expansion of the primary branch}

Write
\begin{equation}\label{branchExpansion}
\begin{split}
w_\varepsilon(x)&=\varepsilon\cos x+a_2\varepsilon^2\cos(2x)+O(\varepsilon^3),\\
\mu(\varepsilon)&=\mu_1^*+\mu_2\varepsilon^2+O(\varepsilon^4), 
\end{split}
\end{equation}
where the remainder is understood in the regularity class supplied by the local bifurcation theorem and the amplitude parameter is chosen so that the coefficient of $\cos x$   {is exactly } $\varepsilon$  {; in particular, the higher-order remainder has no $\cos x$ component. The half-period shift sends the local primary branch with parameter $\varepsilon$ to the same local branch with parameter $-\varepsilon$. Local uniqueness of the Crandall--Rabinowitz branch therefore gives $\mu(-\varepsilon)=\mu(\varepsilon)$, which explains the } absence of a linear term in $\mu$ .\\
For convenience, set
\begin{equation}
t:=\tanh(kh),\qquad
C_1:=\coth(kh)=\frac{1}t,\qquad
C_2:=2\coth(2kh)=\frac{1+t^2}{t}.
\end{equation}
At order $\varepsilon$ equation \eqref{defF} gives the bifurcation condition
\begin{equation}
\mu_1^*C_1=\frac{1}k,
\end{equation}
hence $\mu_1^*=t/k$. Using $T(\cos(nx))=n\coth(nhk)\cos(nx)$, the nonlinear terms needed up to the fundamental harmonic at cubic order are
\begin{equation}\label{nonlinearExpansion}
\begin{split}
w_\varepsilon T w_\varepsilon-[w_\varepsilon T w_\varepsilon]_{2\pi}
&=\frac{C_1}{2}\varepsilon^2\cos(2x)
+\frac{a_2}{2}(C_1+C_2)\varepsilon^3\cos x
+\text{higher harmonics}+O(\varepsilon^4),\\
\mathcal{C}_{kh}(w_\varepsilon w_\varepsilon')
&=\frac{C_2}{4}\varepsilon^2\cos(2x)
+\frac{a_2C_1}{2}\varepsilon^3\cos x
+\text{higher harmonics}+O(\varepsilon^4).
\end{split}
\end{equation}
Here ``higher harmonics'' denotes terms orthogonal to $\cos x$ at cubic order. They are immaterial as they do not contribute to the solvability condition.\\
At order $\varepsilon^2$, the $\cos(2x)$ coefficient gives
\begin{equation}\label{a2eq}
a_2\left(\mu_1^*C_2-\frac{1}k\right)
=\frac{C_1}{2}+\frac{C_2}{4}.
\end{equation}
Since
\begin{equation}
\mu_1^*C_2-\frac{1}k=\frac{t^2}{k},
\end{equation}
we obtain
\begin{equation}\label{a2}
a_2=\frac{k(3+t^2)}{4t^3}.
\end{equation}
At order $\varepsilon^3$, the solvability condition in the fundamental harmonic is
\begin{equation}\label{mu2eq}
\mu_2C_1=\frac{a_2}{2}(2C_1+C_2).
\end{equation}
Using \eqref{a2} and the expressions for $C_1,C_2$, we find
\begin{equation}\label{mu2}
\mu_2=\frac{k(3+t^2)^2}{8t^3}>0
\qquad\text{for every }h,k>0.
\end{equation}
Thus the fixed-$h$ primary branch bends toward larger values of $\mu$ for every finite depth.

\subsection{A negative direction of the Hessian}

The following identity isolates the structural relation between the direction of bifurcation and the Hessian.
\begin{lemma}\label{lemTangentHessian}
Let $\varepsilon\mapsto(\mu(\varepsilon),w_\varepsilon)$ be a differentiable fixed-$h$ solution branch of $\mathcal{F}(\mu,w)=0$. Then
\begin{equation}\label{tangentHessian}
\frac{\delta^2\Lambda}{\delta w^2}(w_\varepsilon,h)
(\dot w_\varepsilon,\dot w_\varepsilon)
=-2g\,\dot\mu(\varepsilon)
\langle\dot w_\varepsilon,Tw_\varepsilon\rangle.
\end{equation}
\end{lemma}
\begin{proof}
Differentiating $\mathcal{F}(\mu(\varepsilon),w_\varepsilon)=0$ and using $\mathcal{F}_\mu(\mu,w)=Tw$ gives
\begin{equation}\label{differentiateBranch}
\mathcal{F}_w(\mu(\varepsilon),w_\varepsilon)\dot w_\varepsilon
+\dot\mu(\varepsilon)Tw_\varepsilon=0.
\end{equation}
Taking the $L^2$ inner product with $\dot w_\varepsilon$ and applying \eqref{HessF} proves \eqref{tangentHessian}.
\end{proof}
From \eqref{branchExpansion} and \eqref{mu2},
\begin{equation}
\dot\mu(\varepsilon)=2\mu_2\varepsilon+O(\varepsilon^3)
\end{equation}
and
\begin{equation}
\langle\dot w_\varepsilon,Tw_\varepsilon\rangle
=\pi\coth(kh)\varepsilon+O(\varepsilon^2),
\end{equation}
and therefore
\begin{equation}\label{negativeTangent}
\frac{\delta^2\Lambda}{\delta w^2}(w_\varepsilon,h)
(\dot w_\varepsilon,\dot w_\varepsilon)
=-4\pi g\,\mu_2\coth(kh)\,\varepsilon^2
+O(|\varepsilon|^3)<0
\end{equation}
for all sufficiently small $\varepsilon\neq0$. Since $\dot w_\varepsilon$ is even and mean-zero, \eqref{negativeTangent} proves variational instability already in the even fixed-$h$ perturbation class, and hence also on the full mean-zero periodic space.

\begin{remark}\label{translationmode}
The full periodic Hessian also has an exact neutral mode generated by translation. Indeed, the functional is invariant under $x\mapsto x+\theta$, and therefore
\begin{equation}
\mathcal{F}_w(\mu(\varepsilon),w_\varepsilon)w_\varepsilon'=0.
\end{equation}
Near the bifurcation this neutral eigenfunction originates from the $\sin x$ direction. This observation is an important consistency check for any second-order spectral perturbation calculation: one of the two eigenvalues emanating from the two-dimensional kernel at the bifurcation point must remain identically zero along the non-trivial branch. The negative direction in \eqref{negativeTangent} is the even amplitude direction and is distinct from this translation mode.
\end{remark}

\begin{remark}
The conclusion is specific to the fixed-conformal-depth Hessian. It does not contradict the shallow-water formal stability result of \cite[Section 4]{constantin2007stability}, because that paper uses a different height-function functional and a different set of constraints, as discussed in the Introduction and Remark \ref{remconstraints}. In particular, the difference $d-h$ in \eqref{relwd} is of order $\varepsilon^2$ along the branch, precisely the order at which the exchange of stability is decided.
\end{remark}

\section{Conclusion and Further Perspectives}\label{SecConclusion}

We have proved variational instability of small-amplitude irrotational traveling waves along the primary bifurcation branch when the conformal mean depth is held fixed. The proof combines the local bifurcation theory with an explicit Lyapunov--Schmidt expansion and an exact identity for the Hessian evaluated on the solution branch. The crucial branch coefficient has been identified as
\begin{equation}
\mu_2=\frac{k(3+\tanh^2(kh))^2}{8\tanh^3(kh)}>0,
\end{equation}
which forces the tangent direction to have negative second variation for sufficiently small non-zero amplitude.\\
For fixed $h$, the trivial solution is strictly variationally stable for $\mu>\tanh(kh)/k$, degenerate at the first bifurcation point, and unstable for smaller $\mu$. At higher bifurcation points the real dimension of the negative subspace is $2(n-1)$, with an additional two-dimensional kernel in the $n$th harmonic.\\
The finite-depth Plotnikov transform  in Section \ref{SecPlotnikov}   {gives } a useful representation of the second variation  {, but unlike its infinite-depth counterpart it contains an explicit finite-rank zero-mode correction. Moreover, } a first-order expansion of the Plotnikov potential is not sufficient to determine a second-order eigenvalue correction when the first-order correction vanishes. The present proof avoids this issue by working directly with the variational equation. Translation invariance also shows that the odd mode generated by $w_\varepsilon'$ remains exactly neutral along the non-trivial branch, whereas the even amplitude mode becomes negative.\\
Finally, the variational notion used here is constraint-dependent and should not be interpreted as a statement about dynamical stability for the Cauchy problem. This distinction is also essential in comparing the result with \cite{constantin2007stability}: their height-function formal stability problem fixes different physical quantities, and physical and conformal mean depths differ at quadratic order along a non-trivial branch. It remains an interesting problem to compare the Morse indices of these different constrained formulations systematically and to determine how the picture changes for larger amplitude or in the presence of vorticity.\\

\textbf{Acknowledgments.} We would like to thank Adrian Constantin for useful discussions and feedback.

\bibliographystyle{abbrv}
\bibliography{Waterwaves}

\end{document}